\documentclass[10pt]{article}
\usepackage{amssymb}
\usepackage{mathrsfs}
\usepackage{CJK}
\usepackage{mathrsfs}
\usepackage{bbm}
\usepackage{amsfonts}
\usepackage{mathrsfs}
\usepackage{bbm}
\usepackage{amscd}
\usepackage{setspace}
\usepackage[top=1in, bottom=1in, left=1.25in, right=1.25in]{geometry}
\makeatletter 
\@addtoreset{equation}{section}
\makeatother  

\newcommand{\mod}{\mbox{\rm mod }}

\newtheorem{thm}{Theorem}[section]

\newtheorem{lem}[thm]{Lemma}
\newtheorem{cor}[thm]{Corollary}

\begin{document}

\baselineskip 20pt
\renewcommand{\baselinestretch}{1.0}

\newcommand{\lr}[1]{\langle #1 \rangle}
\newcommand{\llr}[1]{\langle\hspace{-2.5pt}\langle #1
\rangle\hspace{-2.5pt}\rangle}

\title{ Modulus $p^2$ congruences involving    harmonic numbers }
\author{{Jizhen Yang$^{1}$ and Yunpeng Wang$^{2}$}\\
{\small 1 Department of Mathematics,  Luoyang normal College, }\\
{\small  Luoyang 471934, P. R. China}\\
{\small 2 Department of Mathematics and Physical, Luoyang Institute
of Science and }\\
{\small Technology,}\\
{\small  Luoyang 471023, P. R. China}\\
{\small email:yangjizhen116@163.com} }

\date{}
\maketitle

{\bf Abstract:} The harmonic number
$H_k=\sum_{j=1}^k1/j(k=1,2,3\cdots)$ play an important role in
mathematics. Let $p>3$ be a prime. In this paper, we establish a
number of congruences with the form $\sum_{k=1}^{p-1}k^mH_k^n(\mod
p^2)$ for $m=1,2\cdots,p-2$ and $n=1,2,3$.



\vspace{4mm}

\section{Introduction}

 For $\alpha\in \mathbb{N}$, the
generalized harmonic numbers are  defined by
\begin{eqnarray*}
H_0^{(\alpha)}=0\quad {\rm and}\ \
H_n^{(\alpha)}=\sum_{i=1}^n\frac{1}{i^{\alpha}},\ \ {\rm for}\
n\in\mathbb{N},
\end{eqnarray*}
when $\alpha=1$, they reduce to the well-known harmonic numbers
\begin{eqnarray*}
H_0=0\quad {\rm and}\ \ H_n=\sum_{i=1}^n\frac{1}{i},\ \ {\rm for}\
n\in\mathbb{N}.
\end{eqnarray*}

In 1862, Wolstenholme\cite{W1} proved that if $p>3$ is a prime, then
\begin{eqnarray}\label{Lj8}
H_{p-1}\equiv 0\ \ (\mod p^2)\ \ \ and\ \ \ H_{p-1}^{(2)}\equiv 0\ \
(\mod p).
\end{eqnarray}
 There are many results in related to  the Wolstenholme
theorem(see, e.g., \cite{A}, \cite{A1}, \cite{B}, \cite{S3},
\cite{S4}, \cite{S5}). Recently Z.-W. Sun \cite{S3} gave
\begin{eqnarray*}
&&H_{p-k}\equiv H_{k-1}\ \ (\mod p)\label{L11}
\end{eqnarray*}
and
\begin{eqnarray*}
 &&(-1)^k{p-1\choose k}\equiv
1-pH_k+\frac{p^2}{2}\left(H_k^2-H_k^{(2)}\right)\ \ (\mod
p^3)\label{L12}
\end{eqnarray*}
for every $k=1,2,\ldots p-1$. Using these congruences, Z.-W. Sun
obtained a series of modulus $p$ congruences involving harmonic
numbers. For example,
\begin{eqnarray*}
\sum_{k=1}^{p-1}k^2H_k^2\equiv -\frac{4}{9}\ \ (\mod p)\ \ \ and \ \
\ \sum_{k=1}^{p-1}H_k^2\equiv 2p-2\ \ (\mod p^2).
\end{eqnarray*}

 In this paper, we
investigate the properties of $\sum_{k=1}^{p-1}k^mH_k^n(\mod
p^2)(m=1,2\cdots,p-2,\ n=1,2,3)$
 and establish a series of modulus
$p^2$ congruences involving harmonic numbers, some of which are
generalizations  of partial results in \cite{S3}. For example, we
establish the congruences:
\begin{eqnarray}
&&\sum_{k=1}^{p-1}k^2H_k^2\equiv\frac{79}{108}p-\frac{4}{9}(\mod
p^2) \ \ and\ \ \sum_{k=1}^{p-1}{kH_k^3}\equiv
\frac{27}{8}p-\frac{1}{6}pB_{p-3}-3 (\mod p^2),
\end{eqnarray}
where $B_0,B_1,B_2,\cdots,$ are Bernoulli numbers given by
\begin{eqnarray*}
B_0=1\ \ and\ \ \sum_{k=0}^n{n+1\choose k}B_k=0(n=1,2,3,\cdots).
\end{eqnarray*}

\section{Lemmas}
In this paper, we set
\begin{eqnarray*}
S(m)=\sum_{r=0}^{m}{m+1\choose r}B_rB_{m-r}.
\end{eqnarray*}
 In this section, we
first state some basic facts which will be used very often. For any
prime $p$ and integer $n$, provided $p-1\nmid n$, it is well-known
that
\begin{eqnarray*}
H_{p-1}^{(n)}=\sum_{k=1}^{p-1}\frac{1}{k^n}\equiv0\ \ (\mod p),
\end{eqnarray*}
which implies
\begin{eqnarray*}
H_{p-k}^{(n)}=H_{p-1}^{(n)}-\sum_{j=1}^{k-1}\frac{1}{(p-j)^n}\equiv
(-1)^{n+1}H_{k-1}^{(n)}\ \ (\mod p).
\end{eqnarray*}

\begin{lem}\label{L+1}\rm Let $m, p$ be positive integers. Then
\begin{eqnarray}
\sum_{k=1}^{p-1}k^mH_k&=&\frac{H_{p-1}}{m+1}\sum_{r=0}^m{m+1\choose
r}B_rp^{m+1-r}-(p-1)B_m\nonumber\\
&&-\frac{1}{m+1}\sum_{r=0}^{m-1}\sum_{\lambda=0}^{m-r}\frac{{m+1\choose
r}{m+1-r\choose
\lambda}}{m+1-r}B_rB_{\lambda}p^{m+1-r-\lambda}.\label{l1}
\end{eqnarray}
\end{lem}
{\bf Proof}. For $m, p\in \mathbb{N}$, we have
\begin{eqnarray*}
&&\sum_{k=1}^{p-1}k^mH_k=\sum_{k=1}^{p-1}k^m\sum_{j=1}^k\frac{1}{j}=
\sum_{j=1}^{p-1}\frac{1}{j}\sum_{k=j}^{p-1}k^m=
\sum_{j=1}^{p-1}\frac{1}{j}(\sum_{k=1
}^{p-1}k^m-\sum_{s=1}^{j-1}s^m).
\end{eqnarray*}
It is well-known that
\begin{eqnarray}
&&\sum_{k=1}^{p-1}k^m=\frac{1}{m+1}\sum_{r=0}^{m}{m+1\choose
r}B_rp^{m+1-r} {\rm\ \ \ \ for\  any} \ m, \ p \in \mathbb{Z}^+
.\label{l3}
\end{eqnarray}
(cf. \cite[pp. 230-238]{I}.) Thus
\begin{eqnarray*}
\sum_{k=1}^{p-1}k^mH_k &=&\frac{H_{p-1}}{m+1}\sum_{r=0}^m{m+1\choose
r}B_rp^{m+1-r}- \frac{1}{m+1}\sum_{r=0}^m{m+1\choose
r}B_r\sum_{j=1}^{p-1}j^{m-r}\\
&=&\frac{H_{p-1}}{m+1}\sum_{r=0}^m{m+1\choose
r}B_rp^{m+1-r}-(p-1)B_m- \frac{1}{m+1}\sum_{r=0}^{m-1}{m+1\choose
r}B_r\sum_{j=1}^{p-1}j^{m-r}.
\end{eqnarray*}
Using (\ref{l3}) again, we finally obtain (\ref{l1}).

\begin{lem}\label{L2}\rm Let $m>0$ be an odd integer and $m\neq3$ . Then
\begin{eqnarray}
&&\frac{2}{m+1}\sum_{r=0}^{m}\sum_{\lambda=0}^{m-r}\frac{{m+1\choose
r}{m+1-r\choose
\lambda}}{m+1-r}B_rB_{\lambda}B_{m-r-\lambda}+\frac{1}{m+1}\sum_{r=0}^{m-1}{m+1\choose
r}B_rB_{m-1-r}\nonumber\\
&=&-\frac{m+2}{2m}\sum_{r=0}^{m-1}{m\choose
r}B_{r}B_{m-1-r}+\frac{1}{4}\sum_{r=0}^{m-1}B_rB_{m-1-r}.\label{l5}
\end{eqnarray}
\end{lem}
{\bf Proof}. It is clear for the case of $m=1$. When $m>3$, set
\begin{eqnarray*}
\Sigma&=&\sum_{r=0}^{m}\sum_{\lambda=0}^{m-r}\frac{{m+1\choose
r}{m+1-r\choose \lambda}}{m+1-r}B_rB_{\lambda}B_{m-r-\lambda}\\
&=&\sum_{\begin{array}{l}r=0\\
2\nmid r\end{array}}^{m}\sum_{\lambda=0}^{m-r}\frac{{m+1\choose
r}{m+1-r\choose \lambda}}{m+1-r}B_rB_{\lambda}B_{m-r-\lambda}\\
&&+
\sum_{\begin{array}{l}r=0\\
2\mid r\end{array}}^{m}\sum_{\begin{array}{l}\lambda=0\\
2\nmid \lambda\end{array}}^{m-r}\frac{{m+1\choose r}{m+1-r\choose
\lambda}}{m+1-r}B_rB_{\lambda}B_{m-r-\lambda}\\
&&+\sum_{\begin{array}{l}r=0\\
2\mid r\end{array}}^{m}\sum_{\begin{array}{l}\lambda=0\\
2\mid \lambda\end{array}}^{m-r}\frac{{m+1\choose r}{m+1-r\choose
\lambda}}{m+1-r}B_rB_{\lambda}B_{m-r-\lambda}.
\end{eqnarray*}
Recall that $B_{m}=0$ for $m=3,5,7\cdots$. Then
\begin{eqnarray}
\Sigma&=&\sum_{\lambda=0}^{m-1}\frac{m+1}{m}{m\choose
\lambda}B_1B_{\lambda}B_{m-1-\lambda}+\sum_{\begin{array}{l}r=0\\
2\mid r\end{array}}^{m-1}{m+1\choose r}B_1B_{r}B_{m-1-r}\nonumber\\
&&
+\sum_{\begin{array}{l}r=0\\
2\mid r\end{array}}^{m-1}\frac{m-r}{2}{m+1\choose
r}B_1B_{r}B_{m-1-r}\nonumber.
\end{eqnarray}
When $m>3$ and $m,r$ are odd integers ,we have $B_rB_{m-1-r}=0$.
Thus
\begin{eqnarray}
\Sigma&=&\sum_{r=0}^{m-1}\frac{m+1}{m}{m\choose
r}B_1B_{r}B_{m-1-r}+\sum_{r=0}^{m-1}{m+1\choose r}B_1B_{r}B_{m-1-r}\nonumber\\
&& +\sum_{r=0}^{m-1}\frac{m-r}{2}{m+1\choose
r}B_1B_{r}B_{m-1-r}.\label{l6}
\end{eqnarray}
Therefore, we get
\begin{eqnarray}
&&\frac{2}{m+1}\Sigma+\frac{1}{m+1}\sum_{r=0}^{m-1}{m+1\choose
r}B_rB_{m-1-r}\nonumber\\
&=&-\frac{m+2}{2m}\sum_{r=0}^{m-1}{m\choose
r}B_{r}B_{m-1-r}+\frac{1}{2(m+1)}\sum_{r=0}^{m-1}{m+1\choose
r}B_rB_{m-1-r},\label{l7}
\end{eqnarray}
with the help of (\ref{l6}).
 Observe that
\begin{eqnarray*}
\sum_{r=0}^{m-1}{m+1\choose
r}B_rB_{m-1-r}=\sum_{r=2}^{m-3}{m+1\choose
r}B_rB_{m-1-r}+\frac{m^2+m+2}{2}B_{m-1},
\end{eqnarray*}
since $m$ is an odd integer. Recall that for $n$ is a even
integer(cf. \cite[(4.12)]{D})
\begin{eqnarray*}
(n+2)\sum_{k=2}^{n-2}B_kB_{n-k}=2\sum_{k=2}^{n-2}{n+2\choose
k}B_kB_{n-k}+n(n+1)B_n,
\end{eqnarray*}
which implies
\begin{eqnarray*}
\sum_{r=2}^{m-3}{m+1\choose
r}B_rB_{m-1-r}=\frac{m+1}{2}\sum_{r=2}^{m-3}B_rB_{m-1-r}-\frac{m(m-1)}{2}B_{m-1}.
\end{eqnarray*}
 Hence
\begin{eqnarray}
\sum_{r=0}^{m-1}{m+1\choose r}B_rB_{m-1-r}=
\frac{m+1}{2}\sum_{r=0}^{m-1}B_rB_{m-1-r}\label{l8}.
\end{eqnarray}
In view of  (\ref{l7}) and (\ref{l8}), we finally get
\begin{eqnarray*}
&&\frac{2}{m+1}\Sigma+\frac{1}{m+1}\sum_{r=0}^{m-1}{m+1\choose
r}B_rB_{m-1-r}\nonumber\\
&=&-\frac{m+2}{2m}\sum_{r=0}^{m-1}{m\choose
r}B_{r}B_{m-1-r}+\frac{1}{4}\sum_{r=0}^{m-1}B_rB_{m-1-r}.
\end{eqnarray*}
The proof of (\ref{l5}) is completed.

\begin{lem}\label{L3}\rm Let  $p$ be a prime and $0<m<p-1$ be an even integer. Then
\begin{eqnarray}
\sum_{k=1}^{p-1}\frac{H_k^2}{k^m}\equiv-\sum_{j=0}^{p-1-m}B_jB_{p-1-m-j}
-\sum_{j={p-m}}^{p-2}B_jB_{2p-2-m-j}(\mod p).\label{l9}
\end{eqnarray}
\end{lem}
{\bf Proof}. Observe that
\begin{eqnarray*}
\sum_{k=1}^{p-1}\frac{H_k^2}{k^m}=\sum_{k=1}^{p-1}\frac{H_{p-k}^2}{(p-k)^m}\equiv
\sum_{k=1}^{p-1}\frac{H_{k-1}^2}{k^m}\equiv\sum_{k=1}^{p-1}\frac{(\sum_{j=1}^{k-1}j^{p-2})^2}{k^m}(\mod
p),
\end{eqnarray*}
with the help of Fermat's little theorem. In view of (\ref{l3}), we
obtain
\begin{eqnarray*}
\sum_{k=1}^{p-1}\frac{H_k^2}{k^m}\equiv\sum_{k=1}^{p-1}\frac{1}{k^m}
\left(\frac{1}{p-1}\sum_{j=0}^{p-2}{p-1\choose
j}B_jk^{p-1-j}\right)^2 \equiv\sum_{k=1}^{p-1}\frac{1}{k^m}
\left(\sum_{j=0}^{p-2}B_jk^{p-1-j}\right)^2(\mod p),
\end{eqnarray*}
since
 ${p-1\choose j}B_j\equiv B_j(\mod p)$. Thus,
\begin{eqnarray*}
\sum_{k=1}^{p-1}\frac{H_k^2}{k^m}&\equiv&\sum_{k=1}^{p-1}\frac{1}{k^m}
\left(\sum_{j=0}^{p-2}B_jk^{p-1-j}\right)\left(\sum_{i=0}^{p-2}B_{p-2-i}k^{i+1}\right)\\
&=&\sum_{k=1}^{p-1}\frac{1}{k^m}
\sum_{i,j=0}^{p-2}B_jB_{p-2-i}k^{p-j+i}\\
&=&\sum_{i,j=0}^{p-2}B_jB_{p-2-i}\sum_{k=1}^{p-1}k^{p-j-m+i}(\mod
p).
\end{eqnarray*}
 It is
well-known that for each $t\in\mathbb{Z}$, we have
\begin{eqnarray}
\sum_{k=1}^{p-1}k^t=\left\{\begin{array}{l}-1\ \ (\mod\ p)\ \ {\rm
if}\  p-1\mid
t,\\
\\
0\ \ \ \   (\mod\ p)\ \ {\rm
otherwise}.\\\end{array}\right.\label{l+4}
\end{eqnarray}
(see, e.g., \cite[p.235]{I}.)
Thus, we arrive at
\begin{eqnarray*}
\sum_{k=1}^{p-1}\frac{H_k^2}{k^m}\equiv-\sum_{j=0}^{p-1-m}B_jB_{p-1-m-j}
-\sum_{j={p-m}}^{p-2}B_jB_{2p-2-m-j}(\mod p).
\end{eqnarray*}

\section{Main results}
\begin{thm}\label{T+1}\rm Let  $p>3 $ be a prime and $0<m<p-1$ be a integer. Then
\begin{eqnarray}
&&\sum_{k=1}^{p-1}k^mH_k\equiv B_m-\frac{p}{m+1}S(m)(\mod
p^2).\label{t++1}
\end{eqnarray}
\end{thm}
{\bf Proof}. In view of (\ref{l1}), we have
 \begin{eqnarray*}
\sum_{k=1}^{p-1}k^mH_k\equiv(1-p)B_m
-\frac{1}{m+1}\sum_{r=0}^{m-1}\sum_{\lambda=0}^{m-r}\frac{{m+1\choose
r}{m+1-r\choose \lambda}}{m+1-r}B_rB_{\lambda}p^{m+1-r-\lambda}(\mod
p^2),
\end{eqnarray*}
since $H_{p-1}\equiv0(\mod p^2)$. It is clear that
\begin{eqnarray*}
\sum_{k=1}^{p-1}k^mH_k&\equiv&(1-p)B_m
-\frac{p}{m+1}\sum_{r=0}^{m-1}{m+1\choose r}B_rB_{m-r}\\
&=&B_m-\frac{p}{m+1}S(m)(\mod p^2).
\end{eqnarray*}
The proof of the theorem is completed.
\begin{thm}\label{T1}\rm Let  $p>3 $ be a prime and $0<m<p-1$ be a integer. Then
\begin{eqnarray}
&&\sum_{k=1}^{p-1}k^2H_k^2\equiv\frac{79}{108}p-\frac{4}{9}(\mod
p^2),\label{t+2}\\
&&\sum_{k=1}^{p-1}k^3H_k^2\equiv-\frac{59}{144}p+\frac{1}{6}(\mod
p^2).\label{t1}
\end{eqnarray}
and
\begin{eqnarray}\label{t2}
\sum_{k=1}^{p-1}k^mH_k^2\equiv\left\{\begin{array}{l}B_{m-1}+p\mu_m\
\ (\mod\ p^2)\ \  {\rm if}\ 2\nmid
m {\rm\ and \ }m\neq 3\\
\\
\nu_m+p\lambda_m\ \ \ \ \  \  (\mod\ p^2)\ \ {\rm if}\  2\mid m\
{\rm and}
 \ m\neq2\\ \end{array}\right.,
 \end{eqnarray}
where
\begin{eqnarray*}
&&\mu_m=-\frac{m+2}{2m}S(m-1)+\frac{1}{4}\sum_{r=0}^{m-1}B_rB_{m-1-r},\
\ \ \nu_m=-\frac{2}{m+1}S(m),\\
&&\lambda_m=\frac{2}{m+1}\sum_{r=0}^{m}\frac{{m+1\choose
r}}{m+1-r}B_rS(m-r)-\frac{m^2-m+6}{12}B_{m-2}.
\end{eqnarray*}
\end{thm}
{\bf Proof}. When $m \in \{1,2,\cdots, p-1\}$, we have
\begin{eqnarray*}
\sum_{k=1}^{p-1}k^mH_k^2&=&\sum_{k=1}^{p-1}k^mH_k\sum_{j=1}^k\frac{1}{j}=
\sum_{j=1}^{p-1}\frac{1}{j}\sum_{k=j}^{p-1}k^mH_k
\\ &=&\sum_{j=1}^{p-1}\frac{1}{j}\left(\sum_{k=
1}^{p-1}k^mH_k-\sum_{s=1}^{j-1}s^mH_s\right)\\
&\equiv&-\sum_{j=1}^{p-1}\frac{1}{j}\sum_{s=1}^{j-1}s^mH_s\ (\mod
p^2),
\end{eqnarray*}
since $H_{p-1}\equiv0(\mod p^2)$. In view of (\ref{l1}), we have
\begin{eqnarray*}
\sum_{k=1}^{p-1}k^mH_k^2&\equiv&-\sum_{j=1}^{p-1}\frac{1}{j}\left\{\frac{H_{j-1}}{m+1}\sum_{r=0}^m{m+1\choose
r}B_rj^{m+1-r}-(j-1)B_m\right.\\
&&-\left.
\frac{1}{m+1}\sum_{r=0}^{m-1}\sum_{\lambda=0}^{m-r}\frac{{m+1\choose
r}{m+1-r\choose \lambda}}{m+1-r}B_rB_\lambda
j^{m+1-r-\lambda}\right\}\\
&=&-\frac{1}{m+1}\sum_{r=0}^m{m+1\choose
r}B_r\sum_{j=1}^{p-1}j^{m-r}H_{j-1}+(p-1)B_m-H_{p-1}B_m\\
&&+\frac{1}{m+1}\sum_{r=0}^{m-1}\sum_{\lambda=0}^{m-r}\frac{{m+1\choose
r }{m+1-r\choose
\lambda}}{m+1-r}B_rB_{\lambda}\sum_{j=1}^{p-1}j^{m-r-\lambda}(\mod
p^2).
\end{eqnarray*}
With the help of (\ref{l3}), we can get
\begin{eqnarray*}
\sum_{k=1}^{p-1}k^mH_k^2
&\equiv&-\frac{1}{m+1}\sum_{r=0}^m{m+1\choose
r}B_r\sum_{j=1}^{p-1}j^{m-r}H_{j-1}+\frac{p-1}{m+1}S(m)\\
&&+\frac{p}{m+1}\sum_{r=0}^{m-1}\sum_{\lambda=0}^{m-1-r}\frac{{m+1\choose
r }{m+1-r\choose \lambda}}{m+1-r}B_rB_{\lambda}B_{m-r-\lambda}(\mod
p^2).
\end{eqnarray*}
Note that
\begin{eqnarray*}
\sum_{j=1}^{p-1}j^{m-r}H_{j-1}=\sum_{j=1}^{p-1}j^{m-r}H_{j}-\sum_{j=1}^{p-1}j^{m-1-r}.
\end{eqnarray*}
In view of (\ref{t++1}), we have
\begin{eqnarray*}
\sum_{k=1}^{p-1}k^mH_k^2
&\equiv&-\frac{1}{m+1}\sum_{r=0}^{m}{m+1\choose
r}B_r\left(B_{m-r}-\frac{p}{m+1-r}S(m-r)\right)\\
&&+\frac{p-1}{m+1}S(m)+\frac{p}{m+1}\sum_{r=0}^{m-1}\sum_{\lambda=0}^{m-1-r}\frac{{m+1\choose
r}{m+1-r\choose \lambda}}{m+1-r}B_rB_{\lambda}B_{m-r-\lambda}\\
&&+\frac{1}{m+1}\sum_{r=0}^{m}{m+1\choose
r}B_r\sum_{j=1}^{p-1}j^{m-1-r}\\
&=&-\frac{2}{m+1}S(m)+\frac{2p}{m+1}\sum_{r=0}^{m}\frac{{m+1\choose
r}}{m+1-r}B_rS(m-r)\\
&&+\frac{1}{m+1}\sum_{r=0}^m{m+1\choose
r}B_r\sum_{j=1}^{p-1}j^{m-1-r}(\mod p^2).
\end{eqnarray*}
Observe that
\begin{eqnarray*}
&&\frac{1}{m+1}\sum_{r=0}^m{m+1\choose
r}B_r\sum_{j=1}^{p-1}j^{m-1-r}\\
&=&H_{p-1}B_m+\frac{(p-1)m}{2}B_{m-1}+\frac{1}{m+1}\sum_{r=0}^{m-2}{m+1\choose
r}B_r\sum_{j=1}^{p-1}j^{m-1-r}\\
&\equiv&\frac{(p-1)m}{2}B_{m-1}+\frac{p}{m+1}\sum_{r=0}^{m-2}{m+1\choose
r}B_rB_{m-1-r}\\
&=&-\frac{m}{2}B_{m-1}+\frac{p}{m+1}\sum_{r=0}^{m-1}{m+1\choose
r}B_rB_{m-1-r}(\mod p^2).
\end{eqnarray*}
Therefore
\begin{eqnarray}
\sum_{k=1}^{p-1}k^mH_k^2
&\equiv&\frac{2p}{m+1}\sum_{r=0}^{m}\frac{{m+1\choose
r}}{m+1-r}B_rS(m-r)+\frac{p}{m+1}\sum_{r=0}^{m-1}{m+1\choose
r}B_rB_{m-1-r}\nonumber\\
&&-\frac{2}{m+1}S(m)-\frac{m}{2}B_{m-1}(\mod p^2) \label{t3}.
\end{eqnarray}
Note that $ S(1)=-3/2, S(2)=17/12, S(3)=-5/6$. Taking $m=2, 3$ in
(\ref{t3}), we can obtain  (\ref{t+2}) and (\ref{t1}) respectively.

 When $m>0$ is an odd integer and $m\neq3$ , we have
\begin{eqnarray}
-\frac{2}{m+1}S(m)=-\frac{2}{m+1}\sum_{r=0}^{m}{m+1\choose
r}B_rB_{m-r}=B_{m-1}+\frac{m}{2}B_{m-1}.\label{t4}
\end{eqnarray}
Combining (\ref{l5}), (\ref{t3}) and (\ref{t4}), we obtain
\begin{eqnarray}
\sum_{k=1}^{p-1}k^mH_k^2 \equiv
B_{m-1}-p\left(\frac{m+2}{2m}S(m-1)-\frac{1}{4}\sum_{r=0}^{m-1}B_rB_{m-1-r}\right)(\mod
p^2).\label{t+6}
\end{eqnarray}

 When $m>2$ is an even
integer, we have $B_{m-1}=0$ and
\begin{eqnarray*}
\frac{p}{m+1}\sum_{r=0}^{m-1}{m+1\choose
r}B_rB_{m-r-1}=-\frac{m^2-m+6}{12}pB_{m-2}.
\end{eqnarray*}
Therefore,
\begin{eqnarray*}
\sum_{k=1}^{p-1}k^mH_k^2\equiv \nu_m+p\lambda_m(\mod\ p^2),
\end{eqnarray*}
where
\begin{eqnarray*}
&&\nu_m=-\frac{2}{m+1}S(m),\\
&&\lambda_m=\frac{2}{m+1}\sum_{r=0}^{m}\frac{{m+1\choose
r}}{m+1-r}B_rS(m-r)-\frac{m^2-m+6}{12}B_{m-2}.
\end{eqnarray*}
We are done.

{\bf Remark 1} The congruences (\ref{t+2}) and (\ref{t1}) are the
generalized forms of \cite[(1.2)]{S3} and the third congruences of
\cite[Corollary 1.1]{S3} respectively.

Taking $m=1$ in (\ref{t+6}), we have
\begin{cor} \label{C+1}\rm  Let  $p>3 $ be a prime . Then
\begin{eqnarray}
\sum_{k=1}^{p-1}kH_k^2\equiv -\frac{5}{4}p+1(\mod p^2).\label{cc+2}
\end{eqnarray}
\end{cor}

{\bf Remark 2} The congruences (\ref{cc+2}) is the generalized form
of the first congruences of \cite[Corollary 1.1]{S3}.

Noticing $S(4)=7/90$ and taking $m=5$ in (\ref{t+6}), we have
\begin{cor}\label{C1}\rm Let  $p>5 $ be a prime . Then
\begin{eqnarray}
\sum_{k=1}^{p-1}k^5H_k^2\equiv -\frac{77}{1200}p- \frac{1}{30}(\mod
p^2).\label{t+5}
\end{eqnarray}
\end{cor}

\begin{cor}\label{C2}\rm Let  $p>3 $ be a prime and $0<m<p-1$ be an odd integer. Then
\begin{eqnarray*}
\sum_{k=1}^{p-1}k^mH_k^2\equiv B_{m-1} (\mod p).
\end{eqnarray*}
\end{cor}
{\bf Proof}. In view of (\ref{t3}), we have
\begin{eqnarray*}
\sum_{j=1}^{p-1}k^mH_k^2\equiv-\frac{2}{m+1}\sum_{r=0}^m{m+1\choose
r}B_rB_{m-r}-\frac{m}{2}B_{m-1}=B_{m-1}(\mod p),
\end{eqnarray*}
since $m$ is an odd integer.

{\bf Remark 3} Taking $m=1$ and $m=3$, we also can get the first and
third congruences of \cite[Corollary 1.1]{S3} respectively.

\begin{cor}\label{C3}\rm Let  $p>3 $ be a prime and $0<m<p-1$ be an odd integer. Then
\begin{eqnarray}
\sum_{k=1}^{p-1}\frac{H_k^2}{k^m}\equiv B_{p-2-m} (\mod
p).\label{cc+3}
\end{eqnarray}
\end{cor}
{\bf Proof}. Taking $m\rightarrow p-1-m$ in  Corollary \ref{C2}, we
have
\begin{eqnarray*}
\sum_{k=1}^{p-1}k^{p-1-m}H_k^2\equiv B_{p-2-m} (\mod p),
\end{eqnarray*}
by  Fermat's little theorem, we obtain
\begin{eqnarray*}
\sum_{k=1}^{p-1}\frac{H_k^2}{k^m}\equiv B_{p-2-m} (\mod p)
\end{eqnarray*}
as desired.

Taking $m=4$ in (\ref{t2}), we obtain
\begin{cor}\label{C4}\rm Let  $p>5 $ be a prime. Then
\begin{eqnarray}
\sum_{k=1}^{p-1}k^4H_k^2\equiv \frac{5743}{27000}p-\frac{7}{225}
(\mod p^2).\label{t+4}
\end{eqnarray}
\end{cor}

\begin{cor}\label{C5}\rm Let  $p>5 $ be a prime and $0<m\leq p-5$ is an even integer. Then
\begin{eqnarray}
\sum_{k=1}^{p-1}k^{p-m}H_k^2\equiv B_{p-1-m}-p\left(\frac{m-1}{4}
\sum_{k=1}^{p-1}\frac{H_k^2}{k^m}+\frac{1}{4}\sum_{r=p-m}^{p-2}B_rB_{2p-2-m-r}\right)(\mod
p^2).\label{t5}
\end{eqnarray}
\end{cor}
{\bf Proof}. In view of (\ref{t2}), we have
\begin{eqnarray*}
\sum_{k=1}^{p-1}k^{p-m}H_k^2&\equiv&
B_{p-1-m}+p\left(-\frac{p-m+2}{2(p-m)}S(p-1-m)+\frac{1}{4}\sum_{r=0}^{p-1-m}B_rB_{p-1-m-r}\right)\\
&\equiv&B_{p-1-m}+p\left(\frac{-m+2}{2m}S(p-1-m)+\frac{1}{4}\sum_{r=0}^{p-1-m}B_rB_{p-1-m-r}\right)(\mod
p^2)
\end{eqnarray*}
and
\begin{eqnarray*}
\sum_{k=1}^{p-1}\frac{H_k^2}{k^m}\equiv\sum_{k=1}^{p-1}k^{p-1-m}H_k^2\equiv
-\frac{2}{p-m}S(p-1-m)\equiv\frac{2}{m}S(p-1-m)(\mod p).
\end{eqnarray*}
Combining (\ref{l9})  and the above, we can obtain
\begin{eqnarray*}
\sum_{k=1}^{p-1}k^{p-m}H_k^2\equiv B_{p-1-m}-p\left(\frac{m-1}{4}
\sum_{k=1}^{p-1}\frac{H_k^2}{k^m}+\frac{1}{4}\sum_{r=p-m}^{p-2}B_rB_{2p-2-m-r}\right)(\mod
p^2),
\end{eqnarray*}
we are done.

Taking $m=2$ in  (\ref{t5}), we have
\begin{cor}\label{C6}\rm Let  $p>5 $ be a prime . Then
\begin{eqnarray}
\sum_{k=1}^{p-1}k^{p-2}H_k^2\equiv B_{p-3}(\mod p^2).\label{t6}
\end{eqnarray}
\end{cor}
{\bf Proof}. Observe  that
\begin{eqnarray*}
\sum_{k=1}^{p-1}\frac{H_k^2}{k^2}\equiv 0(\mod p).
\end{eqnarray*}
(cf. \cite[(1.5)]{S3}.) Thus, we obtain (\ref{t6}) from (\ref{t5}).

Taking $m=4$ in (\ref{t5}), we have
\begin{cor}\label{C7}\rm Let  $p>7 $ be a prime . Then
\begin{eqnarray}
\sum_{k=1}^{p-1}k^{p-4}H_k^2\equiv B_{p-5}(\mod p^2).\label{t7}
\end{eqnarray}
\end{cor}
{\bf Proof}. Observe that
\begin{eqnarray*}
\sum_{1\leq i\leq j\leq k \leq p-1}\frac{1}{ij^4k}\equiv
\frac{1}{3}B_{p-3}^2(\mod p),(cf. \cite[Theorem 7.2]{H})
\end{eqnarray*}
 and
\begin{eqnarray*}
\sum_{1\leq i\leq j\leq k \leq
p-1}\frac{1}{ij^4k}&=&\sum_{j=1}^{p-1}\frac{H_j(H_{p-1}-H_{j-1})}{j^4}\equiv
-\sum_{j=1}^{p-1}\frac{H_jH_{j-1}}{j^4}\\
&\equiv&-\sum_{j=1}^{p-1}\frac{H_j^2}{j^4}+\sum_{j=1}^{p-1}\frac{H_j}{j^5}(\mod
p).
\end{eqnarray*}
With the help of  Fermat's little theorem and (\ref{t++1}), we get
\begin{eqnarray*}
\sum_{k=1}^{p-1}\frac{H_k}{k^5}\equiv\sum_{k=1}^{p-1}k^{p-6}H_k\equiv
B_{p-6}=0(\mod p).
\end{eqnarray*}
Combining (\ref{t5}) and the above, we finally obtain
\begin{eqnarray*}
\sum_{k=1}^{p-1}k^{p-4}H_k^2\equiv B_{p-5}(\mod p^2).\
\end{eqnarray*}
The proof of  Corollary \ref{C7} is completed.

\section{Further results }
It is not difficult to  find the following recurrence relation to
$\sum_{k=1}^{p-1}k^{m}H_k^3(\mod p^2)$.

For $m=1,2,\cdots,p-2$, set $a_m\equiv\sum_{k=1}^{p-1}k^{m}H_k^3\
(\mod p^2)$. Then
\begin{eqnarray}
a_{m-1}\equiv\frac{1}{m}\left(-\sum_{i=2}^{m}{m\choose i}a_{m-i}-
3\sum_{k=1}^{p-1}k^{m-1}H_k^2+3\sum_{k=1}^{p-1}k^{m-2}H_k-\sum_{k=1}^{p-1}k^{m-3}\right)(\mod
p^2),\label{r1}
\end{eqnarray}
which can be deduced by
\begin{eqnarray}
\sum_{k=1}^{p-1}(k+1)^mH_k^n\equiv\sum_{k=1}^{p-1}k^mH_{k-1}^n(\mod
p^2).(n\in\mathbb{N})\label{r2}
\end{eqnarray}
Using the Theorem \ref{T+1},  Theorem \ref{T1} and  (\ref{r1}), we
can calculate the congruences $\sum_{k=1}^{p-1}k^{m}H_k^3\ (\mod
p^2)$ for $m=0,1,2,\cdots p-2$. As the examples, we show the cases
of $m=0,1,2,3$.

\begin{thm}\label{T2}\rm Let  $p>3 $ be a prime.  Then
\begin{eqnarray}
&&\sum_{k=1}^{p-1}H_{k}^3\equiv \frac{1}{3}pB_{p-3}-6p+6(\mod
p^2),\label{l19}\\
&&\sum_{k=1}^{p-1}{kH_k^3}\equiv \frac{27}{8}p-\frac{1}{6}pB_{p-3}-3
(\mod p^2),\label{t8}\\
&&\sum_{k=1}^{p-1}k^2H_k^3\equiv-\frac{365}{216}p+\frac{1}{18}pB_{p-3}+\frac{23}{18}
(\mod p^2),\label{c+2}\\
 && \sum_{k=1}^{p-1}{k^3H_k^3}\equiv
\frac{425}{576}p-\frac{5}{12} (\mod p^2)\label{t9}.
\end{eqnarray}
\end{thm}
{\bf Proof of (\ref{l19})}. Taking $m=1$ in (\ref{r1}), we obtain
\begin{eqnarray*}
\sum_{k=1}^{p-1}H_k^3\equiv-
3\sum_{k=1}^{p-1}H_k^2+3\sum_{k=1}^{p-1}\frac{H_k}{k}-H_{p-1}^{(2)}(\mod
p^2).
\end{eqnarray*}
In view of (\ref{r2}), we can get
\begin{eqnarray}
\sum_{k=1}^{p-1}H_k^2\equiv- 2\sum_{k=1}^{p-1}H_k(\mod p^2)\ \ \
and\ \ \ \sum_{k=1}^{p-1}H_k\equiv1-p\ (\mod p^2).\label{cc+1}
\end{eqnarray}
It is easy to find that
\begin{eqnarray*}
H_{p-1}^2=2\sum_{k=1}^{p-1}\frac{H_k}{k}-H_{p-1}^{(2)}.
\end{eqnarray*}
Thus,
\begin{eqnarray*}
\sum_{k=1}^{p-1}\frac{H_k}{k}=\frac{1}{2}\left(H_{p-1}^2+H_{p-1}^{(2)}\right)\equiv\frac{1}{2}H_{p-1}^{(2)}(\mod
p^2).
\end{eqnarray*}
Combining the above, we obtain
\begin{eqnarray*}
\sum_{k=1}^{p-1}H_k^3\equiv\frac{1}{2}H_{p-1}^{(2)}-6p+6\equiv\frac{1}{3}pB_{p-3}-6p+6(\mod
p^2),
\end{eqnarray*}
since $H_{p-1}^{(2)}\equiv2/3pB_{p-3}(\mod p^2)$ by
(\cite[Corollaries 5.1]{S1}).

{\bf Proof of (\ref{t8})}. Taking $m=2$ in (\ref{r1}), we obtain
\begin{eqnarray*}
&&\sum_{k=1}^{p-1}{kH_k^3}\equiv\frac{1}{2}\left(-\sum_{k=1}^{p-1}{H_k^3}-
3\sum_{k=1}^{p-1}{kH_k^2}+3\sum_{k=1}^{p-1}{H_k}-H_{p-1}\right)(\mod
p^2).
\end{eqnarray*}
In view of Corollary {\ref{C+1}} and (\ref{l19}), we have
\begin{eqnarray*}
\sum_{k=1}^{p-1}{kH_k^3}\equiv\frac{27}{8}p -\frac{1}{6}pB_{p-3}-3 \
(\mod p^2),
\end{eqnarray*}
since $\sum_{k=1}^{p-1}H_k\equiv1-p\ (\mod p^2)$.

{\bf Remark 4} The congruences (\ref{t8})   is the generalized form
of the second congruences of \cite[Corollary 1.1]{S3}.

{\bf Proof of (\ref{c+2})}. Taking $m=3$ in (\ref{r1}), we obtain
\begin{eqnarray*}
&&\sum_{k=1}^{p-1}{k^2H_k^3}\equiv\frac{1}{3}\left(-3\sum_{k=1}^{p-1}{kH_k^3}-\sum_{k=1}^{p-1}{H_k^3}-
3\sum_{k=1}^{p-1}{k^2H_k^2}+3\sum_{k=1}^{p-1}{kH_k}-p+1\right)(\mod
p^2).
\end{eqnarray*}
In view of (\ref{t++1}), we can obtain
\begin{eqnarray*}
&&\sum_{k=1}^{p-1}kH_k\equiv\frac{3}{4}p-\frac{1}{2}(\mod p^2).
\end{eqnarray*}
Hence
\begin{eqnarray*}
\sum_{k=1}^{p-1}k^2H_k^3\equiv-\frac{365}{216}p+\frac{1}{18}pB_{p-3}+\frac{23}{18}
\ (\mod p^2),
\end{eqnarray*}
with the help  of (\ref{t+2}), (\ref{l19}) and (\ref{t8}).

{\bf Proof of (\ref{t9})}.  Taking $m=4$ in (\ref{r1}), we obtain
\begin{eqnarray*}
&&\sum_{k=1}^{p-1}{k^3H_k^3}\equiv\frac{1}{4}\left(-6\sum_{k=1}^{p-1}{k^2H_k^3}-4\sum_{k=1}^{p-1}{kH_k^3}-\sum_{k=1}^{p-1}{H_k^3}
-3\sum_{k=1}^{p-1}{k^3H_k^2}+3\sum_{k=1}^{p-1}{k^2H_k}-\sum_{k=1}^{p-1}{k}\right)(\mod
p^2).
\end{eqnarray*}
In view of (\ref{t++1}), we can obtain
\begin{eqnarray*}
&&\sum_{k=1}^{p-1}k^2H_k\equiv-\frac{17}{36}p+\frac{1}{6}(\mod p^2).
\end{eqnarray*}
Hence
\begin{eqnarray*}
 \sum_{k=1}^{p-1}{k^3H_k^3}\equiv
\frac{425}{576}p-\frac{5}{12} (\mod p^2),
\end{eqnarray*}
with the help  of (\ref{t1}), (\ref{l19}), (\ref{t8}) and
(\ref{c+2}).

\end{document}